\input amstex
\documentstyle{amsppt}
\NoBlackBoxes

\hoffset 1.5cm
\voffset1cm
\let\sc\tensc
\define\supp{\operatorname{supp}}
\def\dom{\mathop{\fam0 dom}\nolimits}

\topmatter


\font\bbold=bbold10

\def\reduce{\mskip-5mu }
\def\Bnorml{\reduce\left\bracevert\reduce\vphantom{X}}
\def\Bnormr{\vphantom{X}\reduce\right\bracevert\reduce}

\def\downwardarrow{
\mathord{
 \hbox to 5pt{
\hss$\vcenter{
\hrule\hbox to 2.4pt{
\hss$\mathchar"223$\hss
}
}\hss$
}
}}

\def\tv@rt{{\vert\mkern-2mu\vert\mkern-2mu\vert}}
 \def\tvert#1\tvert{\mathopen\tv@rt#1\mathclose\tv@rt}

\font\ss=lass1000

\thanks
The authors were supported by the Program of Basic Scientific Research of the Siberian Branch
of the Russian Academy of Sciences (Grant No. ~I.1.2, Project 0314--2019--0005).
\endthanks

\title       Two Applications of Boolean Valued Analysis\endtitle

\date August 7, 2019\enddate

\author{A.~G.~Kusraev and S.S. Kutateladze}\endauthor
\address
Southern Mathematical Institute\hfill\break
Vladikavkaz Science Center of the Russian Academy of Sciences\hfill\break
North Ossetian State University, Vladikavkaz, Russia
\email kusraev\@smath.ru\endemail
\endaddress

\address Sobolev Institute of Mathematics\hfill\break
Novosibirsk State University\hfill\break
Novosibirsk, Russia
\email kusraev\@smath.ru;  sskut\@math.nsc.ru\endemail
\endaddress


\abstract
The paper contains two main results that are obtained by Boolean valued analysis.
The first asserts that a universally complete vector lattice
without locally one-dimensional bands can be decomposed into a direct sum of
two vector sublattices that are laterally complete and invariant under all
band projections and there exists a band preserving linear isomorphism of each of these
sublattices to the original lattice. The second result establishes a counterpart of the Ando
Theorem on the joint characterization of~ $A\!L^p$  and ~$c_0$
for the class of cyclic Banach lattices, using the Boolean valued
transfer for injective Banach lattices.
\endabstract

\dedicatory To Yu.~G.~Reshetnyak on the occasion of his 90th birthday
\enddedicatory
\endtopmatter

\document

\head \S\,1.~Introduction \endhead

We present the Boolean valued approach to the two problems of the theory of~vector and
Banach lattices. Prerequisites 
from Boolean valued analysis are given in~Section~2;
details 
can be found in~[1,\,2]. The theory of~vector and Banach lattices is presented
in~[3,\,4].

In~[5, Problem~5], Abramovich and Kitover asked whether
vector lattices~$X$ and $Y$ are lattice isomorphic if there is a~linear invertible
operator $T:X\to Y$ such that $T$ and $T^{-1}$ preserve disjointness. The negative
answer is given in~the same article (see [5, Theorem~13.4]). Some strengthening of~this
result is presented in~Section~3.

The famous Ando Theorem states that a~Banach lattice~$X$ of~dimension~$\geq3$ is
geometrically and lattice isomorphic to~$L^p(\Omega,\Sigma,\mu)$ for~some
$1\leq p\in{\Bbb R}$ and measure space $(\Omega,\Sigma,\mu)$ or to~$c_0(\Gamma)$
for~some nonempty set~$\Gamma$ if and only if each closed sublattice of~$X$
is the image of~a~ contracting positive projection (see, for~example, [4, Theorem~2.7.13]
or~[6, Theorem~1.b.8]). The analogous result in~the class of~${\Bbb B}$-cyclic Banach
lattices is established in~Section~6. The preparatory material about Boolean valued
$AL^p$-spaces and Boolean valued Banach lattices of~the form~$c_0(\Gamma)$
is given in~Sections~4 and~5 respectively.

In~what follows, we let ${\Bbb B}$ stand for a~complete Boolean
algebra with~unity~$\text{\bbold 1}$, zero~$\text{\bbold 0}$, join~$\vee$,
meet~$\wedge$, and complement~$(\cdot)^\ast$, moreover,
$\text{\bbold 1}\ne\text{\bbold 0}$; while ~${\Bbb P}(X)$ is the Boolean algebra of~all
order projections in~a~vector lattice~$X$. By~a~{\it partition of~unity\/} in~${\Bbb B}$
we mean a~family $(b_{\xi})_{\xi\in \Xi}\subset{\Bbb B}$ such that
 $\bigvee\nolimits_{\xi\in \Xi} b_{\xi}=\text{\bbold 1}$ and
$b_{\xi}\wedge b_{\eta}=\text{\bbold 0}$ for $\xi\ne \eta$. The symbol~$:=$ is used
for~assignation by~definition, and~${\Bbb N}$ and~${\Bbb R}$ are the naturals and the
reals respectively.

\head \S\,2.~Preliminaries\endhead

Applying the Transfer and Maximum Principles to~the ${\roman {ZFC}}$-theorem
on~the existence of the reals, find ${\Cal   R}\in{\Bbb V}^{({\Bbb B})}$, called the {\it reals within\/}
${\Bbb V}^{({\Bbb B})}$, satisfying $[\![\,{\Cal   R}$ is the reals$\,]\!]=\text{\bbold 1}$. The following result by Gordon
states that the interpretation of the reals in~the model~${\Bbb V}^{({\Bbb B})}$
is a~universally complete vector~ lattice.

 \proclaim{Theorem 2.1 \rm [7]}
Let ${\Cal   R}$ be the reals within~${\Bbb V}^{({\Bbb B})}$.
Then ${\Cal   R}{\downarrow}$ $($with the descended operations and order$)$
is a~universally complete vector lattice. Moreover, there is a~Boolean isomorphism
$\chi:{\Bbb B}\to {\Bbb P}({\Cal   R}{\downarrow})$
such that
$$
 \chi(b)x=\chi(b)y\ \Longleftrightarrow\ b\leq [\![x=y]\!],\quad
 \chi(b)x\leq\chi(b)y\ \Longleftrightarrow\ b\leq [\![x\leq y]\!]
                                                                \tag1
$$
for~all $x,y\in{\Cal   R}{\downarrow}$ and $b\in{\Bbb B}$.
\endproclaim

 {\sc Proof.}  See~[3, Theorem~8.1.2; 2, Theorem~2.2.4].\qed

\medskip
Consider the reals ${\Bbb R}\in{\Bbb V}$. The standard name
 ${\Bbb R}^{\scriptscriptstyle\wedge}$ is a~field within
 ${\Bbb V}^{({\Bbb B})}$, and we may assume that
 $[\![{\Bbb R}^{\scriptscriptstyle\wedge}$ is a~dense subfield
of~${\Cal   R}]\!]=\text{\bbold 1}$ (see [2, Subsections~2.2.2 and~ 2.2.3]).
In~some questions, it is important to~know when
${\Bbb R}^{\scriptscriptstyle\wedge}={\Cal   R}$
(see, for~example, [8]). The following result by Gutman
answers this question in terms of~equivalent algebraic properties~${\Bbb B}$
and ${\Cal   R}{\downarrow}$:

 \proclaim{Theorem~2.2 \rm [9]}
Let ${\Bbb B}$ be a~complete Boolean algebra and let ${\Cal   R}$ be the
reals within~${\Bbb V}^{({\Bbb B})}$. Then the following are equivalent:

$(1)$~${{\Bbb V}}^{({\Bbb B})}\!\models\!{\Cal   R}
={{\Bbb R}}^{\scriptscriptstyle\wedge}$.

 $(2)$~${\Bbb B}$ is $\sigma$-distributive;

 $(3)$~${\Cal   R}{\downarrow}$ is locally one-dimensional.
\endproclaim

 {\sc Proof.} See [3, Theorem~5.1.6; 2, Theorems~4.4.9 and~4.7.6].
 \qed

Recall the notion of ${\Bbb B}$-cyclic Banach lattice (see~[1,\,2]).

\demo{Definition 2.3}
An~order projection~$\pi$ in~a~Banach lattice~$X$ is
an~$M$-{\it projection\/} if $\|x\|=\max\{\|\pi x\|,\allowmathbreak \|\pi^\ast x\|\}$ for all
$x\in X$, where $\pi^\ast\!:=I_X-\pi$. The set of~all $M$-projections in~$X$
is denoted by~${\Bbb M}(X)$. A ~{\it complete Boolean algebra of~$M$-projections\/} in~$X$
is a~subset ${\Cal   B}\subset{\Bbb M}(X)$ that is an~order closed subalgebra
in~the complete Boolean algebra~${\Bbb P}(X)$.
\enddemo

Note that ~${\Bbb M}(X)$ is a~subalgebra of~${\Bbb P}(X)$ in~$X$ but can fail to~be order complete
(see [10, Theorem~1.10]).

 \demo{Definition 2.4}
If $(b_\xi)_{\xi \in \Xi}$ is a~partition of~unity in~${\Cal   B}$ and
$(x_\xi)_{\xi \in \Xi}$ is a~family in~$X$ then $x\in X$ satisfying
$b_\xi x_\xi=b_\xi x$ for all $\xi \in\Xi$ is called
the {\it mixing\/} of~$(x_\xi)$ with  ~$(b_\xi)$. 
A~Banach lattice~$X$
is called ${\Bbb B}$-{\it cyclic\/} if ${\Bbb B}$ is the complete Boolean algebra
isomorphic to a~complete Boolean algebra~${\Cal   B}$ of~$M$-projections
in~$X$ and the mixing of every family in~the unit ball of~$X$ with~respect to every
partition of~unity in~${\Cal   B}$ (with~the same index set) exists and belongs to the unit ball.
\enddemo

In~what follows, we identify~${\Bbb B}$ and ${\Cal   B}$ and assume that
 ${\Bbb B}\subset{\Bbb P}(X)$ (see [11, Definition 2.5; 2, Definition~5.7.13]).
We say that ${\Bbb B}$-cyclic Banach lattices~$X$ and $Y$ are {\it ${\Bbb B}$-isometric\/}
and write $X\simeq_{{\Bbb B}}Y$ if there exists a~lattice isometric isomorphism
between~$X$ and $Y$ commuting with the elements of~${\Bbb B}$.

Thus, a~Banach lattice is ${\Bbb B}$-cyclic if it is a~${\Bbb B}$-cyclic Banach space
with~respect to~some Boolean algebra of~$M$-projections
${\Bbb B}\subset{\Bbb P}(X)$
(cf. [3, Definitions~ 7.3.1 and~ 7.3.3] and [2, Definition~ 5.8.8]).

 \demo{Definition 2.5}
Denote by $\Lambda={\Cal   R}{\Downarrow}$ the bounded part of~the universally
complete vector lattice~ ${\Cal   R}{\downarrow}$, i.e.,~ $\Lambda$ is the order dense
ideal in~${\Cal   R}{\downarrow}$ generated by~the order unity
 ${\text{\bbold 1}}\!:=1^{\scriptscriptstyle\wedge}\in{\Cal   R}\!\!\downarrow$.
Take a~nonzero Banach space
${\Cal   X}\!=({\Cal   X},\|\cdot\|_{\Cal   X})$
within ${\Bbb V}^{({\Bbb B})}$ and put
 $$
 \gathered
 {\Cal   X}{\Downarrow}:=\{x\in{\Cal   X}{\downarrow}:
 \Bnorml x\Bnormr\in\Lambda\},\quad \tvert x\tvert\!:=\|\Bnorml
 x\Bnormr\|_\infty\quad(x\in X),
 \\
 \|\lambda\|_\infty:=\inf\{\alpha>0: |\lambda|\leq\alpha\text{\bbold 1}\},
\quad \lambda \in\Lambda ,
 \endgathered
 $$
 where $\Bnorml\cdot\Bnormr$ is the descent of~the norm $\|\cdot\|_{\Cal   X}$,
i.e.,
 $[\![\Bnorml x\Bnormr=\|x\|_{\Cal   X}]\!]=\text{\bbold 1}$ for~all
 $(x\in{\Cal   X}{\downarrow})$.
Then
${\Cal   X}{\Downarrow}\!:=  ({\Cal   X}{\Downarrow},\tvert\cdot\tvert)$
is called the {\it bounded descent\/} of~ ${\Cal   X}$.
\enddemo

\medskip
Since $\Lambda$ is an~order complete $A\!M$-space with~unity,
${\Cal   X}{\Downarrow}$ is a~Banach space with~mixed norm over~$\Lambda$ and so
a~${\Bbb B}$-cyclic Banach space (see [3, 7.3.3]). The following result states that
the concept of ${\Bbb B}$-cyclic Banach lattice is nothing but interpretation
of~the notion of Banach lattice in~a~Boolean valued model.

\proclaim{Theorem 2.6 \rm [11]}~A~Banach lattice~$X$ is ${\Bbb B}$-cyclic if and only
if ~${\Bbb V}^{({\Bbb B})}$ contains a~Banach lattice~${\Cal   X}$ unique
up~to~a~lattice isometry whose bounded descent is ${\Bbb B}$-isometric to~$X$.
Moreover, $\pi\mapsto\pi{\Downarrow}\!:=\pi{\downarrow}|_X$ is
an~isomorphism of the Boolean algebras ${\Bbb M}({\Cal   X}){\downarrow}$ and~${\Bbb M}(X)$
{\rm (}in~symbols:
 ${\Bbb M}({\Cal   X}){\downarrow}\simeq{\Bbb M}({\Cal   X}{\Downarrow}))$.
 \endproclaim

 {\sc Proof.} See [11, Theorem~2.1] or [2, Theorem~5.9.1]. \qed

 \demo{Definition 2.7}
The Banach lattice~${\Cal   X}$ within~${\Bbb V}^{({\Bbb B})}$ of~Theorem~2.6
is called the {\it Boolean valued realization\/} of~a~${\Bbb B}$-cyclic Banach lattice~$X$.
\enddemo

\demo{Remark 2.8}
The bounded descent of~Definition~2.5 appeared firstly under another name
in the articles~[12,\,13] by Takeuti in his study of  von~Neumann algebras and
$C^\ast$-algebras by using Boolean valued models (also see~[14]).
\enddemo

\medskip
We will need the following result on~the structure of~Boolean valued cardinals:

\proclaim{Proposition 2.9}
$x\in{\Bbb V}^{({\Bbb B})}$ is a~cardinal within~${\Bbb V}^{(\Bbb B)}$
if and only if there are a~nonempty set of~cardinals $\Gamma\in{\Bbb V}$ and
a~partition of~unity $(b_\gamma)_{\gamma\in\Gamma}\subset {\Bbb B}$ such that
 $x=\mathop{\roman{mix}}
\nolimits_{\gamma\in\Gamma} b_\gamma\gamma^{\scriptscriptstyle\wedge}$ and
$\gamma^{\scriptscriptstyle\wedge}$ is a~cardinal within~${\Bbb V}^{(\Bbb B_\gamma)}$,
where ${\Bbb B}_\gamma\!:=[\text{\bbold 0},\,b_\gamma]$ and $b_\gamma\ne0$ for all
$\gamma\in\Gamma$.
\endproclaim

 {\sc Proof.} See ~[1, Theorem~9.1.3]; also  [2,~1.9.11]. \qed

 \head\S\,3.~Contracting Linear Isomorphisms\endhead

 In~[5, Section~6, 13], we propose an~approach to~constructing some counterexamples
to~the above-mentioned Problem~B that uses the method of~$d$-basis. In~this
section, we show that, up to~passage to~a~suitable Boolean-valued model, this
is equivalent to~the application of~a~classical Hamel basis.

 \demo{Definition 3.1}
Suppose that $X$ is a~vector lattice and $u\in X$. An~element $v\in X$
is called {\it {\rm a}~fragment\/} or  {\it {\rm a}~component\/} of~$v$ if $|v|\wedge|u-v|=0$.
The set of~all fragments of~$u$ is denoted by~${\Bbb C}(u)$. A~subset~$X_0$ is
called {\it fragment closed\/} if ~${\Bbb C}(u)$ lies in~$X_0$ for all $u\in X_0$ (see~[5, Proposition~ 4.9]).
\enddemo

 \demo{Definition 3.2}
A~sublattice $X_0\subset X$ is called {\it laterally complete\/} if each disjoint set
of~positive elements in~$X_0$ has a~supremum; and ${\Bbb P}(X)$-{\it invariant}, if
$\pi(X_0)\subset X_0$ for~all $\pi\in{\Bbb P}(X)$. An~operator $T:X_0\to  X$ is called
{\it contracting\/} or {\it band preserving\/} if $T(B\cap X_0)\subset B$ for every
band $B\subset X$ (see~[2, Definition~ 4.1.2]).
\enddemo

\proclaim{Lemma 3.3}
Let ${\Cal   R}$ be the reals within ${\Bbb V}^{({\Bbb B})}$. Consider
the universally complete vector lattice $X\!:={\Cal   R}{\downarrow}$. The following
are equivalent for~a~sublattice $X_0\subset X$:

 {\rm(1)}~$X_0$ is order complete, laterally complete, and fragment closed;

 {\rm(2)}~$X_0$ is order complete, laterally complete, and ${\Bbb P}(X)$-invariant;

 {\rm(3)}~$X_0$ is laterally complete, ${\Bbb P}(X)$-invariant, and
$X_0^{\perp\perp}=X$;

 {\rm(4)}~$X_0={\Cal   X}_0{\downarrow}$ for~some vector sublattice
 ${\Cal   X}_0$ of~the field ${\Cal   R}$ regarded as a~vector lattice
 over the subfield~${\Bbb R}^{\scriptscriptstyle\wedge}$.
 \endproclaim

 {\sc Proof.} This is immediate from~[2, Theorem~2.5.1]. \qed

 \proclaim{Lemma 3.4}
Let ${\Bbb P}$ be a~proper subfield of~${\Bbb R}$. Then there exist ${\Bbb P}$-linear
subspaces~${\Cal   X}_1$ and ${\Cal   X}_2$ in~${\Bbb R}$ such that ${\Cal   X}_k$
and ${\Bbb R}$ are isomorphic as vector spaces over~${\Bbb P}$ but not isomorphic
as ordered vector spaces over~${\Bbb P}$. Moreover,
${\Bbb R}={\Cal   X}_1\oplus{\Cal   X}_2$ and the corresponding projections
 $p_1:{\Bbb R}\to{\Cal   X}_1$ and $p_2:{\Bbb R}\to{\Cal   X}_2$ are not order
bounded.
 \endproclaim

 \demo{Proof}
The real~${\Bbb R}$ is a~finite extension of~no~proper subfield
${\Bbb P}\subset{\Bbb R}$ (see, for~example, [15, Lemma~17]). Consequently,
${\Bbb R}$ is an~infinite-dimensional space over~${\Bbb P}$. Let ~${\Cal   E}$
be a~Hamel basis for the ${\Bbb P}$-vector space~${\Bbb R}$ and let~ $|{\Cal   E}|$
be the cardinality of~${\Cal   E}$. Since the cardinal~$|{\Cal   E}|$ is infinite,
we can choose a~proper subset ${\Cal   E}_1\varsubsetneq{\Cal   E}$ so that~
${\Cal   E}_1$ and  ${\Cal   E}_2\!:={\Cal   E}\setminus{\Cal   E}_1$ have the same
cardinality equal to~$|{\Cal   E}|$, i.e., $|{\Cal   E}|=|{\Cal   E}_1|=|{\Cal   E}_2|$.
Let ${\Cal   X}_k$ denote the ${\Bbb P}$-linear subspace in~${\Bbb R}$ generated
by~${\Cal   E}_k$, where $k=1,2$. Then
 $\{0\}\varsubsetneq{\Cal   X}_k\varsubsetneq{\Bbb R}$, where
${\Cal   X}_k$ and ${\Bbb R}$ are isomorphic as vector spaces over~${\Bbb P}$ since
$|{\Cal   E}|=|{\Cal   E}_k|$. If ~${\Cal   X}_k$ and~ ${\Bbb R}$
were isomorphic as ordered vector spaces over~${\Bbb P}$ then ~${\Cal   X}_k$ would
be order complete and so we would obtain the contradictory equality
${\Cal   X}_k={\Bbb R}$. By the choice of~${\Cal   E}_1$ and~
${\Cal   E}_2$, we have ${\Bbb R}={\Cal   X}_1\oplus{\Cal   X}_2$, and so there are
projections $p_k:{\Bbb R}\to{\Cal   X}_k$, where $p_1+p_2=I_{\Cal   R}$. If
$p_k:{\Bbb R}\to{\Cal   X}_k$ is order bounded then ~$p_k$ is continuous as an~additive
function in~${\Bbb R}$; consequently, $p_k(x)=c_kx$ $(x\in{\Bbb R})$ for~some
$c_k\in{\Bbb R}$. But then $c^2_k=1$; i.e., either $c_k=1$ and then
${\Cal   X}_k={\Bbb R}$ or $c_k=0$ and then ${\Cal   X}_k=\{0\}$; in~both cases, we get
a~contradiction. \qed
\enddemo

\smallskip
Let us now prove the first main result of~the article which is the interpretation
of Lemma~3.4 in~an~arbitrary Boolean valued model.

\proclaim{Theorem~3.5}
Let $X$ be a~universally complete lattice not containing nonzero locally one-dimensional
bands. Then there are fragment closed laterally complete vector sublattices
$X_1\subset X$ and $X_2\subset X$ and linear bijections $T_1:X_1\to X$ and $T_2:X_2\to X$ such that

 {\rm (1)}~$X=X_1\oplus X_2$ and $X=X_1^{\perp\perp}=X_2^{\perp\perp}$;

 {\rm (2)}~$T_k$ and $T_k^{-1}$ preserve bands $(k=1,2)$;

 {\rm (3)}~the canonical projections $\pi_1:X\to X_1$ and $\pi_2:X\to X_2$ preserve
bands;

 {\rm (4)}~none of~the sublattices~$X_1$ and $X_2$ is order complete and so is not lattice
isomorphic to~$X$.
\endproclaim

 \demo{Proof}
By~the Gordon Theorem~2.1, we may assume without loss of~generality that
$X={\Cal   R}{\downarrow}$. Since~ $X$ contains no~locally one-dimensional bands,
by the Gutman Theorem~2.2,
 $[\![{\Cal   R}\ne{\Bbb R}^{\scriptscriptstyle\wedge}]\!]=\text{\bbold 1}$.
The Transfer Principle enables us to~apply Lemma~3.4 within~${\Bbb V}^{({\Bbb B})}$;
therefore, there are ${\Bbb R}^{\scriptscriptstyle\wedge}$-linear
subspaces~${\Cal   X}_1$ and ${\Cal   X}_2$ in~${\Cal   R}$ such that
 ${\Cal   R}={\Cal   X}_1\oplus{\Cal   X}_2$ as well as
 ${\Bbb R}^{\scriptscriptstyle\wedge}$-linear isomorphisms~$\tau_k$
in~${\Cal   X}_k$ on~${\Cal   R}$; moreover, ${\Cal   X}_k$ and ${\Cal   R}$
are nonisomorphic as ordered vector spaces
over~${\Bbb R}^{\scriptscriptstyle\wedge}$ $(k=1,2)$. Let
 $p_k:{\Cal   R}\to{\Cal   X}_k$ $(k=1,2)$ be the projections corresponding
to the decomposition ${\Cal   R}={\Cal   X}_1\oplus{\Cal   X}_2$. Put
$X_k\!:={\Cal   X}_k{\downarrow}$, $T_k\!:=\tau{\downarrow}$,
 $S_k\!:=\tau_k^{-1}{\downarrow}$, and $P_k\!:=p_k{\downarrow}$. Since
 ${\Cal   X}_k$ is linearly isomorphic to~${\Cal   R}$,
 $[\![{\Cal   X}_k\ne\{0\}]\!]=\text{\bbold 1}$. Hence, $X_k^{\perp\perp}=X$.
By~Lemma~3.3, $X_1$ and $X_2$ are fragment closed and laterally
complete. Moreover, $X=X_1\oplus X_2$. The operators~$S_k$,~ $T_k$, and~ $P_k$ preserve
bands and are ${\Bbb R}$-linear by~[2, Theorem~4.3.4]. Furthermore,
 $S_k=(\tau_k{\downarrow})^{-1}=T_k^{-1}$ and~ $P_1$ and~ $P_2$ are the projections
corresponding to~the decomposition $X=X_1\oplus X_2$. It remains to~observe that
$X_k$ and $X$ are lattice isomorphic if and only if ${\Cal   X}_k$ and
 ${\Cal   R}$ are isomorphic as ordered vector spaces
over~${\Bbb R}^{\scriptscriptstyle\wedge}$ and $P_1$ and $P_2$ are order bounded
if and only if $p_1$ and $p_1$ are order bounded within~${\Bbb V}^{({\Bbb B})}$. \qed
\enddemo

 \demo{Remark 3.6}
Lemma 3.4 and hence Theorem 3.5 are based on~the absorption property of~infinite
cardinals: $\varkappa+\varkappa=\varkappa$. However, this property holds also
in~the case when the number of summands is infinite but does not exceed~$\varkappa$;
i.e., $\varkappa=\sum\nolimits_{\alpha\in{\roman A}}\varkappa_\alpha$
if $\varkappa_\alpha=\varkappa$ for all $\alpha\in{\roman A}$ and
 $|{\roman A}|\leq\varkappa$. Consequently, we arrive at the version of~Theorem~3.5
with infinitely many sublattices.
\enddemo

\head \S\, 4.~Boolean Valued $AL^p$-Spaces\endhead

In~this section, we introduce the class of~${\Bbb B}$-cyclic Banach lattices
that are classical $AL^p$-spaces in~the Boolean valued model. Start from the key
notion of injective Banach lattice.

 \demo{Definition 4.1}
A~real Banach lattice~$X$ is called {\it injective\/} if for every  Banach lattice~$Y$,
every closed vector sublattice $Y_0\subset Y$, and every  positive linear operator
$T_0:Y_0\to X$, there is a~positive linear extension $T:Y\to X$ of~$T_0$ such
that $\|T_0\|=\|T\|$. In~other words, the injective Banach lattices are injective objects
in~the category of~Banach lattices with~contracting positive operators as~morphisms
(see~ [2, \S\,5.10; 11]).
\enddemo

\medskip
Recall some facts on~the structure of~Banach lattices. If $X\ne\{0\}$ is an~injective
Banach lattice then~ ${\Bbb M}(X)$ is an~order closed subalgebra in~${\Bbb P}(X)$
(see Definition~2.3). In~this case, $X$ is a~${\Bbb B}$-cyclic Banach lattice for~any
order closed subalgebra ${\Bbb B}\subset{\Bbb M}(X)$. This circumstance makes it
possible to~establish a~Boolean valued Transfer Principle for~injective Banach lattices.

 \proclaim{Lemma 4.2}
For every injective Banach lattice~$X$, there is
 ${\Cal   X}\in{\Bbb V}^{({\Bbb B})}$, where ${\Bbb B}={\Bbb M}(X)$, such that
$[\![{\Cal   X}$ is an~$AL$-space\,$]\!]=\text{\bbold 1}$ and $X$ is isometrically and
lattice ${\Bbb B}$-isomorphic to~the bounded descent ${\Cal   X}{\Downarrow}$.
\endproclaim

 {\sc Proof.} See~[11, Theorem 4.4]. \qed

 \demo{Definition 4.3}
A~positive operator~$T$ from a~vector lattice~$X$ into a~vector lattice~$Y$ is said
to~possess the {\it Levy property\/} whenever $Y=T(X)^{\perp\perp}$ and
$\sup x_\alpha$ exists in~$X$ for an~ increasing chain $(x_\alpha)\subset X_+$
if only $(Tx_\alpha)$ is order bounded in~$Y$. A~ {\it Maharam operator\/}
is an~interval preserving order continuous operator; i.e., $T([0,x])=[0,Tx]$ for all
$x\in X_+$ (see [3, 3.4.1]). We say that $T$ is {\it strictly positive\/} if~ $T$
is positive and $T(|x|)=0$ implies $x=0$.

\smallskip
Denote by $(\Lambda({\Bbb B}),\|\cdot\|_\infty)$ an~order complete $A\!M$-space
with~unity~$\text{\bbold 1}$ such that the Boolean algebras~${\Bbb B}$ and
 ${\Bbb P}(\Lambda({\Bbb B}))$ are isomorphic. The norm~$\|\cdot\|_\infty$ is introduced
in~the same way as in~Definition~2.5:
$\|x\|_\infty=\inf\{\alpha>0:|x|\leq\alpha\text{\bbold 1}\}$.
\enddemo

 \proclaim{Lemma 4.4}
For~an~arbitrary injective Banach lattice $(X,\|\cdot\|)$ with Boolean algebra
of~$M$-projections ${\Bbb B}={\Bbb M}(X)$, there exists a~strictly positive
Maharam operator $\Phi:X\to\Lambda({\Bbb B})$ with~the Levy property such that
$\|x\|=\|\Phi(|x|)\|_\infty$ for~all $x\in X$. Moreover, there exists a~Boolean
isomorphism~$h$ from~${\Bbb P}(\Lambda({\Bbb B}))$ onto~${\Bbb M}(X)$ such that
 $\pi\circ\Phi=\Phi\circ h(\pi)$ for~all $\pi\in{\Bbb P}(\Lambda({\Bbb B}))$.
\endproclaim

 {\sc Proof.} See~[11, Corollary~4.5]. \qed

\smallskip
If $X$ and $\Phi$ are the same as in~Lemma~4.4 then we will write $X=L^1(\Phi)$.
Assume that $X=L^1(\Phi)$, ${\Bbb B}={\Bbb M}(X)$, and $\Lambda\!:=\Lambda({\Bbb B})$. Let
$X^u$ and $\Lambda^u$ be the universal completions of~$X$ and~ $\Lambda$.
Take $\text{\bbold 1}\leq p\in\Lambda^u$ and put
$$
\gathered
 X^p\!:=L^p(\Phi)\!:=\{x\in X^u: |x|^p\in X,\
 \Bnorml x\Bnormr_p\!:=(\Phi(|x|^p))^{\frac1{p}}\in\Lambda\},\\
 \|x\|_p\!:=\|\Bnorml x\Bnormr_p \|_\infty\quad(x\in X^p).
\endgathered
$$
Note that if $\text{\bbold 1}\leq p\in\Lambda^u$ then $p^{-1}$ exists and belongs
to~$\Lambda$. Moreover, we may assume that $\Lambda^u$ is a~sublattice in~$X^u$ with
the same order unity as~$X^u$. Then $X^u$ is an~$f$-algebra and $\Lambda^u$ is an
$f$-subalgebra of~$X^u$; therefore, the mappings $x\mapsto|x|^p$ and
$\lambda\mapsto|\lambda|^{1/p}$ are well defined on~$X^u$ and $\Lambda^u$.

 \proclaim{Lemma 4.5}
Let $X$ be an~injective Banach lattice and let ${\Cal   X}$ be the representation of~ $X$
in~the Boolean valued model~${\Bbb V}^{({\Bbb B})}$, where
 ${\Bbb B}={\Bbb M}(X)$, and $\text{\bbold 1}\leq p\in\Lambda^u$. Then
$[\![{\Cal   X}^p$ is an~$AL^p$-space$]\!]=\text{\bbold 1}$ and
$X^p$ is an~${\Bbb B}$-cyclic Banach lattice, where ${\Cal   X}^p{\Downarrow}$ and
$X^p$ are isometrically and lattice ${\Bbb B}$-isomorphic.
\endproclaim

 \demo{Proof}
By~the Gordon Theorem~2.1 and Lemma~4.2, we may assume that
 $\Lambda^u={\Cal   R}{\downarrow}$ and $X={\Cal   X}{\Downarrow}$. By~Lemma~4.4,
there exists a~strictly positive Maharam operator $\Phi:X\to\Lambda({\Bbb B})$
with~the Levy property such that $\|x\|=\|\Phi(|x|)\|_\infty$ for~all $x\in X$. Put
$\varphi\!:=\Phi{\uparrow}$ and note that $[\![\,\varphi:{\Cal   X}\to{\Cal   R}$
is an~order continuous strictly positive functional with~the Levy property
$]\!]=\text{\bbold 1}$; moreover, $X$ is an~order ideal in~${\Cal   X}{\downarrow}$
and the restriction of~$\varphi{\downarrow}$ to~${\Cal   X}{\Downarrow}$
coincides with~$\Phi$ (see [2, Theorem 5.2.8]). If ${\Cal   X}^u$ is the universal
completion of~${\Cal   X}$ within~${\Bbb V}^{{\Bbb (B)}}$ then the vector
lattices~${\Cal   X}^u{\downarrow}$ and $X^u$ can be identified since they are
lattice isomorphic by~[2, Theorem~2.11.8]. Moreover, the multiplication
in~$X$ coincides with~the descent of~the multiplication in~${\Cal   X}$; hence,
$x\mapsto|x|^p$ $(x\in X)$ is the descent of~the mapping $x\mapsto|x|^p$
 $(x\in{\Cal   X})$. Since $p\in{\Cal   R}{\downarrow}$, by~the Transfer Principle,
$[\![{\Cal   X}^p$ is the $AL^p$-space$]\!]=\text{\bbold 1}$. By~definition, the norm
$\tvert\cdot\tvert_p\in{\Bbb V}^{({\Bbb B})}$ of the $AL^p$-space~ ${\Cal   X}^p$ has
the form $\tvert x\tvert_p=(\varphi(|x|^p))^{1/p}$ $(x\in{\Cal   X}^p)$. Using the rules
for~the descent of~the composition of~operators and the inverse operator (see
 [2, 1.5.5]) as well as the equality $\Phi=\varphi{\downarrow}|_X$, we infer that
$\Bnorml\cdot\Bnormr_p=\tvert\cdot\tvert_p{\downarrow}$. It~ remains to~notice that
the containments $x\in X^p$ and
$\Bnorml x\Bnormr_p\in\Lambda$ are equivalent for~all~ $x\in{\Cal   X}^u{\downarrow}$. \qed
\enddemo

 \proclaim{Theorem 4.6}
For every injective Banach lattice $X=L^1(\Phi)$ and every $1\leq p\in\Lambda^u$,
the space $X^p=L^p(\Phi)$ is a~${\Bbb B}$-cyclic Banach lattice, and the Boolean algebras
${\Bbb B}={\Bbb M}(X^p)$ and ${\Bbb M}(X)$ are isomorphic.
\endproclaim

 {\sc Proof.} This is immediate from~Theorem~2.6 and Lemma~4.5. \qed

 \demo{Remark 4.7}
Applying the Gutman Theorem to the Banach--Kantorovich representation
(see [3, Theorem~2.4.10]), we can show that $L^p(\Phi)$ is isometrically and lattice
isomorphic to~the Banach lattice of~continuous sections of~a~continuous Banach
bundle of~$AL^p$-spaces (hence, of~the classical Lebesgue spaces)
$(L^{p(\omega)}(\varphi_\omega))_{\omega\in\Omega}$, where $\Omega$ is the Stone
compact space of the Boolean algebra~${\Bbb B}$. Details can be found in~[3, \S\,4 and~\S\,5].
A~similar result on~the representation of~injective Banach lattices by~means
of~a~continuous bundle of~$AL$-spaces was obtained by~Haydon in~[16, Theorem~7B].
\enddemo

 \demo{Remark 4.8}
Suppose that $X$ is an~injective Banach lattice,
 ${\Bbb B}={\Bbb M}(X)$, and $\Lambda=\Lambda({\Bbb B})$, while~$Q$ and~ $P$
are the Stone compact spaces of the Boolean algebras~${\Bbb P}(X)$ and ${\Bbb B}$
respectively. The inclusion ${\Bbb M}(X)\subset{\Bbb P}(X)$ induces a~continuous
epimorphism $\tau:Q\to P$. Then there exists a~modular Maharam measure
$\mu:{\Cal   B}\to\Lambda=C(P)$, where ${\Cal   B\!:=}{\Cal   B}(Q)$ is the Borel
$\sigma$-algebra of~$Q$, such that $X$ is lattice and isometrically ${\Bbb B}$-isomorphic
to~$L^1(\mu)\!:=L^1(Q,{\Cal   B},\mu)$ (see [16, Theorem~6H]). Here to~the Maharam
operator~$\Phi$  of~Lemma~4.4 there corresponds the integration operator
 $f\mapsto\int\nolimits_Qf\,d\mu$ $(f\in L^1(\mu))$. We can now define  ~$L^p(\mu)$
as the space of~functions $\mu$-integrable with~a~variable exponent
$p\in C_\infty (P)$:
 $$
 L^p(\mu)=\bigg\{f\in L^0(\mu): \int\limits_Q
 |f(q)|^{p(\tau(q))}\,d\mu(q)\in\Lambda\bigg\}.
 $$
The order metric properties of~classes of~functions integrable with~variable exponent
with~respect to~a~vector deserve extra study.
\enddemo

 \head\S\,5. The Banach Lattices $c_0(\Gamma)$ and~$C_\#(Q,c_0(\Gamma))$\endhead

Recall that $l^\infty(\Gamma)$ stands for~the Banach lattice of~all bounded functions
$x:\Gamma\to{\Bbb R}$ with~the norm
 $\|x\|_\infty\!:=\sup\nolimits_{\gamma\in\Gamma}|x(\gamma)|$, and the Banach
lattice~$c_0(\Gamma)$ is defined as~the closure in~$l^\infty(\Gamma)$
of~the sublattice~$c_{00}(\Gamma)$ consisting of~functions with~finite support
 $\supp(x)\!:=\{\gamma\in\Gamma:x(\gamma)\ne0\}$.

 \proclaim{Lemma 5.1}
Let $\Gamma$ be an~arbitrary nonempty set. Then
 $c_0(\Gamma^{\scriptscriptstyle\wedge})$ is the completion
within~${\Bbb V}^{({\Bbb B})}$ of~the metric space
 $c_0(\Gamma)^{\scriptscriptstyle\wedge}$.
\endproclaim

 \demo{Proof}
Denote by~${\Cal   P}_{\operatorname{fin}}(\Gamma)$ the set of~all finite
subsets of~a~set~$\Gamma$. We will use the formula
$$
 {\Bbb V}^{({\Bbb B})}\models\text{\ss P} _{\operatorname{fin}}
(X^{\scriptscriptstyle\wedge})=
 \text{\ss P} _{\operatorname{fin}} (X)^{\scriptscriptstyle\wedge}
                                                                \tag2
$$
 (see [1, 5.1.9]). Let
$c_{00}(\Gamma^{\scriptscriptstyle\wedge},{\Bbb R}^{\scriptscriptstyle\wedge})$
be the subset of~$c_0(\Gamma^{\scriptscriptstyle\wedge})$ consisting of~functions
 $x:\Gamma^{\scriptscriptstyle\wedge}\to{\Bbb R}^{\scriptscriptstyle\wedge}$
with~finite support $\supp(x)$. Since $c_{00}(\Gamma)$ is dense in~$c_0(\Gamma)$
and
$[\![c_{00}(\Gamma^{\scriptscriptstyle\wedge},{\Bbb R}^
{\scriptscriptstyle\wedge})$
is dense in~$c_0(\Gamma^{\scriptscriptstyle\wedge})]\!]=\text{\bbold 1}$,
it suffices to~show that
$[\![c_{00}(\Gamma^{\scriptscriptstyle\wedge},{\Bbb R}^
{\scriptscriptstyle\wedge})
 \subset c_{00}(\Gamma)^{\scriptscriptstyle\wedge}]\!]=\text{\bbold 1}$.
 By the properties of~descents (see [2, 1.5.2]), the last can be rewritten as
 $$
[\![(\forall x\in c_{00}(\Gamma^{\scriptscriptstyle\wedge},
 {\Bbb R}^{\scriptscriptstyle\wedge}))\,x\in
 c_{00}(\Gamma)^{\scriptscriptstyle\wedge}]\!]=
 \bigwedge\{[\![x\in c_{00}(\Gamma)^{\scriptscriptstyle\wedge}]\!]:
 [\![x\in c_{00}(\Gamma^{\scriptscriptstyle\wedge},
 {\Bbb R}^{\scriptscriptstyle\wedge})]\!]=\text{\bbold 1}\}=\text{\bbold 1}.
 $$
Take $x\in c_{00}(\Gamma^{\scriptscriptstyle\wedge},
 {\Bbb R}^{\scriptscriptstyle\wedge}){\downarrow}$ and consider the modified
descent
$x\downwardarrow:\Gamma\to{\Bbb R}^{\scriptscriptstyle\wedge}{\downarrow}$
 (see [2, 1.5.8]). Verify that
$[\![x\in c_{00}(\Gamma)^{\scriptscriptstyle\wedge}]\!]=\text{\bbold 1}$.
Involving~(2), we infer
$$
\gathered
 [\![x\in c_{00}(\Gamma^{\scriptscriptstyle\wedge},
 {\Bbb R}^{\scriptscriptstyle\wedge})]\!]=
 [\![(\exists\theta\in
\text{\ss  P}_{\operatorname{fin}}(\Gamma^{\scriptscriptstyle\wedge}))\,
 (\forall\gamma\in\Gamma^{\scriptscriptstyle\wedge})\,
 (\gamma\notin\theta\rightarrow x(\gamma)=0)]\!]
 \\
 =\bigvee\limits_{\theta\in\text{\ss P}_{\operatorname{fin}} (\Gamma)}
 [\![(\forall\gamma\in\Gamma^{\scriptscriptstyle\wedge})\,
 (\gamma\notin\theta^{\scriptscriptstyle\wedge}\rightarrow x(\gamma)=0)]\!].
\endgathered
$$
By~the Exhaustion Principle (see [2, 1.2.8]), there exists a~partition of~unity
$(b_\theta)_{\theta\in\Theta}$ of ${\Bbb B}$, where
 $\Theta\!:=\text{\ss P}_{\operatorname{fin}}(\Gamma)$, such that for~all
 $\theta\in\Theta$  we have
$$
\gathered
 b_\theta\leq[\![(\forall\gamma\in\Gamma^{\scriptscriptstyle\wedge})\,
 (\gamma\notin\theta^{\scriptscriptstyle\wedge}\rightarrow x(\gamma)=0)]\!]
\\
=\bigwedge_{\gamma\in\Gamma}
[\![\gamma^{\scriptscriptstyle\wedge}\notin\theta^{\scriptscriptstyle\wedge}]\!]
 \Rightarrow[\![x(\gamma^{\scriptscriptstyle\wedge})=0)]\!]
=\bigwedge_{\gamma\in\Gamma\setminus\theta}[\![x{\downwardarrow}(\gamma)=0)]\!].
\endgathered
$$
Since $\{x{\downwardarrow}(\gamma):\gamma\in\theta\}$ is a~finite subset
of~${\Bbb R}^{\scriptscriptstyle\wedge}\!{\downarrow}$, there are a~partition
of~unity $(\pi_\xi)_{\xi\in\Xi}$ in~${\Bbb B}$ and a~family
$(t_{\theta,\xi})_{\xi\in\Xi}$
such that
$b_\theta x{\downwardarrow}(\gamma)=
 \mathop{\roman{mix}}
\nolimits_{\xi\in\Xi}\pi_\xi t_{\theta,\xi}^{\scriptscriptstyle\wedge}$
for~all $\xi\in\Xi$ and $\gamma\in\theta$. Define  $x_{\xi}\in c_{00}(\Gamma)$,
by~setting $x_{\xi}(\gamma)=t_{\xi,\gamma}$ if $\gamma\in\theta$ and
$x_{\xi}(\gamma)=0$ otherwise. Then
$[\![x_{\xi}^{\scriptscriptstyle\wedge}\in
 c_{00}(\Gamma)^{\scriptscriptstyle\wedge}]\!]=\text{\bbold 1}$ and
 $$
\gathered
[\![x\in c_{00}(\Gamma)^{\scriptscriptstyle\wedge}]\!]
\geq[\![x_{\xi}^{\scriptscriptstyle\wedge}\in
 c_{00}(\Gamma)^{\scriptscriptstyle\wedge}]\!]
 \wedge[\![x_{\xi}^{\scriptscriptstyle\wedge}=x]\!]
=[\![x_{\xi}^{\scriptscriptstyle\wedge}=x]\!]
 \geq\bigwedge_{\gamma\in\Gamma}
 [\![x_{\xi}^{\scriptscriptstyle\wedge}(\gamma^{\scriptscriptstyle\wedge})
 =x(\gamma^{\scriptscriptstyle\wedge})]\!]
 \\
=\bigwedge_{\gamma\in\theta}
 [\![t_{\xi}^{\scriptscriptstyle\wedge}
 =x{\downwardarrow}(\gamma)]\!]
 \wedge\bigwedge_{\gamma\notin\theta}
 [\![0=x{\downwardarrow}(\gamma)]\!]
\geq\pi_{\xi}\wedge b_\theta
\endgathered
$$
for~all $\xi\in\Xi$ and $\theta\in\Theta$;
whence
$[\![x\in c_{00}(\Gamma)^{\scriptscriptstyle\wedge}]\!]\geq
 \bigvee\nolimits_{\xi\in\Xi}b_\theta\wedge\pi_{\xi}=b_\theta$.
Consequently,
 $[\![x\in c_{00}(\Gamma)^{\scriptscriptstyle\wedge}]\!]
 \geq\bigvee\nolimits_{\theta\in\Theta}=\text{\bbold 1}$;  q.e.d. \qed
\enddemo

 \demo{Definition 5.2}
Suppose that $Q$ is an~extremally disconnected compact space and $X$ is a~Banach lattice.
Denote by~$C_\# (Q,X)$ the set of cosets of continuous
vector-functions $u:\dom(u)\subset Q\to X$ such that $Q\setminus\dom(u)$
is a~meager in~$Q$ and the continuous extension $\Bnorml u\Bnormr$
of~the pointwise norm $q\mapsto\|u(q)\|$ to~the whole of~$Q$ belongs to~the Banach
lattice~$C(Q)$ of~continuous functions. Vector-functions~$u$ and $v$ are equivalent
if $u(q)= v(q)$ for~all $q\in\dom(u)\cap\dom(v)$. If $\tilde{u}$ is the coset
of~$u$ then we put
$\Bnorml\tilde{u}\Bnormr\!:=\Bnorml{u}\Bnormr$.
The set~$C_\# (Q,X)$ is naturally endowed with~the structure of~a~module over
the $f$-algebra~$C(Q)$ and the norm $\|u\|\!:=\|\Bnorml u\Bnormr\|_\infty$
(cf.~[3, 2.3.3]).
\enddemo

\proclaim{Lemma 5.3}
Let $X\in {\Bbb V}$ be a~Banach lattice and let ${\Cal   X}$ be~the completion of~the metric
space~$X^{\scriptscriptstyle\wedge}$ within ${\Bbb V}^{({\Bbb B})}$. Then $[\![\,{\Cal   X}$
is the Banach lattice\,$\,]\!]=\text{\bbold 1}$ and ${\Cal   X}{\Downarrow}$ is
${\Bbb B}$-cyclic Banach lattice isometrically and lattice ${\Bbb B}$-isomorphic
to~$C_\#(Q,X)$, where $Q$ is the Stone compact space of~the Boolean algebra~${\Bbb B}$.%
\endproclaim

 \demo{Proof}
We must only apply the description of~the descent ${\Cal X}{\downarrow}$
of~[3, Theorem~8.3.4(1)] and carry out the obvious passage to~the bounded parts. \qed
\enddemo

 \proclaim{Corollary 5.4}
$C_\#(Q,X)$ is a~${\Bbb B}$-cyclic Banach lattice, where ${\Bbb B}$ is isomorphic
to~the Boolean algebra of~clopen subsets in~$Q$. Here the $M$-projection in~$C_\#(Q,X)$
corresponding to~a~clopen set $U\subset Q$ is the multiplication
by~the characteristic function~$\chi_U$.
 \endproclaim

 \proclaim{Corollary 5.5}
Let $\Gamma$ be a~nonempty set and let $Q$ be the Stone compact space of a~complete
Boolean algebra~${\Bbb B}$. Then the ${\Bbb B}$-cyclic Banach lattices
$c_0(\Gamma^{\scriptscriptstyle\wedge}){\Downarrow}$ and $C_\#(Q,c_0(\Gamma))$ are lattice and isometrically ${\Bbb B}$-isomorphic.
\endproclaim

 {\sc Proof.} This is immediate from Lemmas~5.1 and~5.3. \qed

 \demo{Remark 5.6}
Lemma~5.3 is a~particular case of~the general result on~the functional realization
of~the construction of~the Boolean extension which was obtained by Gordon and Lyubetsky (see~[17,\,18]).
\enddemo

 \demo{Remark 5.7}
Important information on~the structure of~a~Banach lattice is given by~the possibility
of embedding into~it of~the classical sequence spaces ~$c_0$, $l^1$, and
$l^\infty$ (see, for~example, [4, \S\,2.3 and~ \S\,2.4]). Corollary~5.5 shows that
the analogous role in~the theory of~${\Bbb B}$-cyclic Banach lattices must belong
to~the ${\Bbb B}$-embeddability of~the Banach lattices~$C_\#(Q,c_0)$, $C_\#(Q,l^1)$,
and $C_\#(Q,l^\infty)$ (see also Remarks~6.5 and~6.6 ~below).
\enddemo

  \head \S\,6.~The Boolean Valued Ando Theorem \endhead

We begin with formulating the Ando Theorem whose Boolean valued interpretation
leads us to~an~analogous result for~${\Bbb B}$-cyclic Banach lattices.

 \proclaim{Theorem 6.1 \rm [19]}
Let $X$ be a~Banach lattice of~dimension~$\geq3$. Then $X$ is isometrically isomorphic
to~$L^p(\Omega,\Sigma,\mu)$ for~some $1\leq p\in{\Bbb R}$ and measure space
$(\Omega,\Sigma,\mu)$ or to~$c_0(\Gamma)$ for~some nonempty set~$\Gamma$ if and only
if each closed vector sublattice in~$X$ is the image of~a~positive contracting
projection.
\endproclaim

Recall that a {\it  contracting positive projection\/} in~a~Banach lattice
is a~positive linear operator $P:X\to X$ with $P^2=P$ and $\|P\|\leq1$.

\demo{Definition 6.2}
A~subspace~$X_0$ in~a~${\Bbb B}$-cyclic Banach lattice~$X$ is called
${\Bbb B}$-{\it complete\/} if $X_0$ contains the mixings of all families
from~$X_0$ by all  partitions of~unity in~${\Bbb B}$. We say that
the Boolean dimension ${\Bbb B}\text{-}\dim (X)$ is at least ~3
if we can choose three elements in~$X$ so that the least closed ${\Bbb B}$-complete
subspace containing them is ${\Bbb B}$-isomorphic to~$\Lambda({\Bbb B})^3$.
\enddemo

 \proclaim{Lemma 6.3}
If ${\Cal   X}\in{\Bbb V}^{({\Bbb B})}$ is the Boolean valued realization
of~a~${\Bbb B}$-cyclic Banach lattice~$X$ then
 $$
 {\Bbb B}\text{-{\rm dim}}(X)\geq3\ \Longleftrightarrow\
 [\![\text{\rm dim}({\Cal   X})\geq3^{\scriptscriptstyle\wedge}]\!]=\text{\bbold 1}.
 $$
 \endproclaim

\demo{Proof}
From [2, Proposition~ 5.8.5, Theorem~ 5.8.11] it follows that a~${\Bbb B}$-cyclic Banach
lattice is an~$f$-module over the ring $\Lambda\!:=\Lambda({\Bbb B})$ (see [2, 2.11.1]
for the definition of $f$-module). Take $e_1,e_2,e_3\in X$ and consider the mapping
$\varphi:(\lambda_1,\lambda_2,\lambda_3)\mapsto\lambda_1e_1+
 \lambda_2e_2+\lambda_3e_3$
from~$\Lambda^3$ into~$X$. It is easy to~see that if~ $e_1$,~ $e_2$, and~$e_3$ are
$\Lambda$-linearly independent then $X_0\!:=\varphi(\Lambda^3)$ is
a~${\Bbb B}$-complete subspace in~$X$. At~the same time, $\varphi$ is injective
if and only if~ $e_1$, $e_2$, and~$e_3$ are $\Lambda$-linearly independent. It remains
to~observe that the $\Lambda$-linear independence of~$e_1$, $e_2$, and~$e_3$ is equivalent
to~the equality
 $[\![\text{\rm
dim}({\Cal   X})\geq3^{\scriptscriptstyle\wedge}]\!]=\text{\bbold 1}$. \qed
\enddemo

\smallskip
We have all ingredients necessary for~formulating and proving
the second main result of this article.

\proclaim{Theorem 6.4}
Let ${\Bbb B}$ be a~complete Boolean algebra and let $Q$ be the Stone compact space of~${\Bbb B}$.
The following are equivalent for~a~${\Bbb B}$-cyclic Banach lattice~$X$ satisfying
the condition ${\Bbb B}\text{-{\rm dim}}(X)\geq3$:

 {\rm(1)}~there is a~ contracting positive projection onto any ${\Bbb B}$-complete
closed sublattice in~$X$ that commutes with the projections from~${\Bbb B}$;

 {\rm(2)}~there is a~partition of~unity $(\pi_\gamma)_{\Gamma\cup\{0\}}$
in~${\Bbb B}$, with ~$\Gamma$ a~nonempty set of~cardinals, such that
 $\pi_0X\simeq_{\pi_0{\Bbb B}}L^p(\Phi)$ for~some
 $\text{\bbold 1}\leq p\in\Lambda({\Bbb B})^u$
and an~injective Banach lattice $L\!:=L^1(\Phi)$ with
${\Bbb M}(L)\simeq\pi_0{\Bbb B}$, and
 $\pi_\gamma X\simeq_{\pi_\gamma{\Bbb B}}C_\#(Q_\gamma,c_0(\gamma))$ for~all
 $\gamma\in\Gamma$, where $Q_\gamma$ is a~clopen subset of~$Q$ corresponding
to~the projection~$\pi_\gamma$.
\endproclaim

 \demo{Proof}
(1)~$\Longrightarrow$~(2): In~correspondence with~Theorem~2.6, we may assume that
$X={\Cal   X}{\Downarrow}$, where~ ${\Cal   X}$ is a~Banach lattice
within~${\Bbb V}^{({\Bbb B})}$; moreover, the relations
${\Bbb B}\text{-{\rm dim}}(X)\geq3$ and
$[\![\dim({\Cal   X})\geq3^{\scriptscriptstyle\wedge}]\!]=\text{\bbold 1}$
are equivalent by~Lemma~6.3. It is easy to~check that ${\Cal   X}_0$ is a~closed
sublattice in~${\Cal   X}$ if and only if $X_0\!:={\Cal   X}_0{\Downarrow}$
is a~${\Bbb B}$-complete norm closed sublattice in~$X$. The operator $P:X\to X$
has the form $P=\pi{\Downarrow}$ for~an~operator $\pi:{\Cal   X}\to{\Cal   X}$ if and
only if $P$ is extensional (see [2, 1.5.6]), and by~(1) the last means that $P$
commutes with~the projections from~${\Bbb B}$. Moreover, $P$ is a~
contracting positive  projection only if $[\![\,\pi$ is a~ contracting positive projection$\,]\!]=\text{\bbold 1}$.
Thus, Proposition~6.4(1) is equivalent to
$[\![\,$each closed sublattice in~${\Cal   X}$ admits a~ contracting positive
projection$\,]\!]=\text{\bbold 1}$.
By the above-mentioned theorem, 6.4(1) is equivalent also to~the fact that
${\Cal   X}$ is lattice and isometrically isomorphic to~$L^p(\mu)$ for~some
$1\leq p\in{\Cal   R}$ and some measure space $(\Omega,\Sigma,\mu)$ or to~$c_0(S)$
for~some nonempty set~$S$ within~${\Bbb V}^{({\Bbb B})}$. Write this down
in~symbolic form:
 $$
 [\![(\exists\,1\leq p\in{\Cal   R})(\exists\,\varphi){\Cal   X}\simeq
 L^p(\varphi)\vee(\exists\,S){\Cal   X}\simeq c_0(S)]\!]=\text{\bbold 1}.
 $$
This implies the existence of~two pairwise complementary elements
 $\pi_0,\pi_0^\ast\allowmathbreak \in{\Bbb B}$ such that if we put
${\Bbb B}_1\!:=[0,\pi_0]$ and ${\Bbb B}_2\!:=[0,\pi^\ast_0]$ then we can choose
$p,\varphi\in{\Bbb V}^{({\Bbb B}_1)}$ and $S\in{\Bbb V}^{({\Bbb B}_2)}$
for which the following two assertions hold:

 $(\ast)$~within the model ${\Bbb V}^{({\Bbb B}_1)}$, the real $1\leq p\in{\Cal   R}$,
the injective Banach lattice~$L$, and the strictly positive order continuous functional
$\varphi:L\to{\Cal   R}$ with~the Levy property satisfy the estimate
 $\pi_0\leq[\![{\Cal   X}\simeq L^p(\varphi)]\!]$;

 $(\ast\ast)$~within the model~${\Bbb V}^{({\Bbb B}_2)}$, the set~$S$
satisfies the estimate $\pi^\ast_0\leq[\![{\Cal   X}\simeq c_0(S) ]\!]$.

Put $\Phi\!:=\varphi{\Downarrow}$ and observe that, by~Lemma~4.5,
 $L^p(\Phi)\simeq_{\pi_0{\Bbb B}}L^p(\varphi){\Downarrow}$. Now, assertion~$(\ast)$ shows
that $\pi_0X\simeq_{\pi_0{\Bbb B}}L^p(\Phi)$.

Let $\varkappa\in{\Bbb V}^{({\Bbb B}_2)}$ be the cardinality of~$S$. In~view
of~Proposition~2.9, the Boolean valued cardinal~$\varkappa$ is a~mixing of standard cardinals; i.e.,
$\varkappa=\mathop{\roman{mix}}\nolimits_{\gamma\in\Gamma}
 b_\gamma\gamma^{\scriptscriptstyle\wedge}$,
where $\Gamma$ is a~nonempty set of~cardinals, $(b_\gamma)_{\gamma\in\Gamma}$
is a~partition of~unity in~${\Bbb B}_2$, $\gamma^{\scriptscriptstyle\wedge}$
is a~cardinal within ${\Bbb V}^{(\Bbb B_\gamma)}$,
${\Bbb B}_\gamma\!:=[\text{\bbold 0},\,b_\gamma]$ for~all $\gamma\in\Gamma$. Passing
to~the model ${\Bbb V}^{(\Bbb B_\gamma)}$ (by~the scheme of~[2, 1.3.7]) and involving
assertion~$(\ast\ast)$, we see that
${\Cal   X}\simeq
c_0(\gamma^{\scriptscriptstyle\wedge})$.
Applying Lemma~5.1 and Definition~5.2, we see that
 $b_\gamma X\simeq_{b_\gamma{\Bbb B}}C_\#(Q_\gamma,c_0(\gamma))$.

\medskip
 (2)~$\Longrightarrow$~(1)~As in~Lemma~3.3, the mapping
 ${\Cal   X}_0\mapsto{\Cal   X}_0{\Downarrow}$ performs a~bijection between
 the norm-closed sublattices of~${\Cal   X}$ and the ${\Bbb B}$-complete the norm-closed sublattices
of~$X$. Similarly, the contracting positive projections in~${\Cal   X}$ are
in~a~one-to-one correspondence with the contracting positive projections in~${\Cal   X}$
commuting with projections in~${\Bbb B}$. We are left with noting that if assertion~(2)
holds then the Banach lattice~${\Cal   X}$ is either an~$AL^p$-space or~$c_0(\Gamma)$,
and so any closed sublattice of~${\Cal   X}$ admits a~contracting positive projection
by~the Transfer Principle. \qed
\enddemo

 \demo{Remark 6.5}
A~nonzero element~$x$ in~a~${\Bbb B}$-cyclic Banach lattice~$X$ is called
a~${\Bbb B}$-{\it atom\/} if for every  pair of~disjoint elements $x,y\in X_+$, $y+z\leq x$,
there exists a~projection $\pi\in{\Bbb B}$ such that $\pi y=0$ and $\pi^\ast z=0$.
The lattice~$X$  is called ${\Bbb B}$-{\it atomic\/} if zero is the only element
in~$X$ disjoint from every ${\Bbb B}$-atom in~$X$. If $X$ is the same as in~Theorem~6.4(2)
then there exists an~$M$-projection $\rho\leq\pi_0$ such that $\rho^\ast X$ has
no~${\Bbb B}$-atoms and $\rho X$ admits the representation
 $$
 \rho X\simeq_{\rho{\Bbb B}}\bigg(\sum\limits_{\gamma
 \in\Delta}^\oplus C_\#(P_\gamma,
l^{p(\gamma)}(\gamma))\bigg)_{l_\infty},
 $$
where $\Delta$ is a~set of~cardinals and $(P_\gamma)_{\gamma\in\Delta}$ is a~family
of~extremally disconnected compact spaces. This claim can be deduced
from [20, Theorem~5.4].
\enddemo

 \demo{Remark 6.6}
If a~${\Bbb B}$-cyclic Banach lattice~$X$ satisfies the conditions of~Lemma~6.4(1)
then, as in~6.4(2), $\pi_0X\simeq_{\pi_0{\Bbb B}}L^p(\Phi)$, and if
 $\pi\!:=\pi_0^\ast=I_X-\pi_0$ then we have the representation
$$
 \pi X\simeq_{\pi{\Bbb B}}\bigg(\sum\limits_{\gamma
 \in\Gamma}^\oplus C_\#(Q_\gamma, c_0(\gamma))\bigg)_{l_\infty}.
 $$
Here the family of~extremal compact spaces $(Q_\gamma)_{\gamma\in\Gamma}$ (as in~Remark~6.5)
is  nonunique in~general. The reason is the phenomenon of~the
 {\it cardinal shift\/} in~a~Boolean valued model. Overcoming
this difficulty and obtaining a~unique representation is possible in~the same way as
in~[20], by involving the notions of pure ${\Bbb B}$-{\it atomicity\/} and {\it stability\/}
[20, Definitions~3.4 and~3.6].
\enddemo

\smallskip
The authors express their gratitude to the referee for~removing numerous typos and inaccuracies.

 \Refs

 {

 \ref\no 1
 \by    Kusraev~A.~G. and Kutateladze~S.~S.
      \book Introduction to Boolean Valued Analysis
      \publ     Nauka
      \publaddr Moscow
      \yr   2005
      \lang Russian
\endref

 \ref\no 2
 \by Kusraev~A.~G. and Kutateladze~S.~S.
 \book Boolean Valued Analysis: Selected Topics {\rm (Trends Sci. South Russia)}
 \publaddr Vladikavkaz
 \publ Vladikavkaz Sci. Center Press
 \yr 2014
 \endref

 \ref\no 3
  \by    Kusraev~A. G.
  \book    Dominated Operators
  \publaddr   Dordrecht
  \publ    Kluwer Academic Publishers
  \yr       2001
\endref

 \ref\no 4
 \by Meyer-Nieberg~P.
 \book Banach Lattices
 \publaddr Berlin etc.
 \publ Springer
 \yr 1991
 \endref

 \ref\no 5
 \by Abramovich~ Y.~ A. and  Kitover A.~K.
 \book Inverses of Disjointness Preserving Operators
\publaddr Providence
\publ Amer. Math. Soc.
  \yr 2000) (Mem. Amer. Math. Soc.; No.~679
 \endref

 \ref\no 6
 \by Lindenstrauss J. and Tzafriri L.
 \book Classical Banach Spaces. V.~2. Function Spaces
 \publaddr Berlin etc.
\publ Springer-Verlag
 \yr 1979
  \endref

 \ref\no 7
  \by Gordon~E.~I.
\paper Real numbers in Boolean-valued models of set theory and $K$-spaces
\jour Dokl. Akad. Nauk SSSR
\yr 1977
\vol 237
\issue 4
\pages 773--775
\endref

 \ref\no 8
 \by Gutman A.~E., Kusraev A.~G., and Kutateladze S.~S.
 \paper The Wickstead problem
 \jour Sib. \`Elektron. Mat. Izv.
 \yr 2008
 \vol 5
 \pages 293--333
 \endref

 \ref\no 9
 \by Gutman A.~E.
 \paper Locally one-dimensional $K$-spaces and $\sigma$-distributive Boolean algebras
 \jour Sib. Adv. Math.
 \yr 1995
 \vol 5
 \issue 2
 \pages 99--121
 \endref

\ref\no 10
 \by Harmand~P., Werner~D., and Wener~W.
 \book $M$-Ideals in Banach Spaces and Banach Algebras
 \publ Springer-Verlag
 \publaddr Berlin etc.
 \yr 1993) (Lecture Notes Math.; V.~1547
 \endref

 \ref\no 11
 \by Kusraev A.~G.
 \paper The Boolean transfer principle for injective Banach lattices
 \jour Sib. Mat. Zh. 
 \yr 2015
 \vol 25
 \issue 1
 \pages 57--65 
 \endref

 \ref\no 12
 \by Takeuti~G.
 \paper Von Neumann algebras and Boolean valued analysis
 \jour J.~Math. Soc. Japan.
 \yr 1983
 \vol 35
 \issue 1
 \pages 1--21
 \endref

 \ref\no 13
 \by Takeuti~G.
 \paper $C^*$-algebras and Boolean valued analysis
 \jour Japan J. Math.
 \yr 1983
 \vol 9
 \issue 2
 \pages 207--246
 \endref

 \ref\no 14
 \by Ozawa~M.
 \paper Boolean valued interpretation of Banach space theory and module structure of von Neumann algebras
 \jour Nagoya Math.~ J.
 \yr 1990
 \vol 117
 \pages 1--36
 \endref

 \ref\no 15
 \by Coppel W.~A.
 \book Foundations of Convex Geometry
 \publaddr Cambridge
 \publ Cambr. Univ. Press
 \yr 1988
 \endref

 \ref\no 16
 \by Haydon R.
\paper Injective Banach lattices
\jour Math. Z.
\vol 156
\pages 19--47
yr 1977
 \endref

 \ref\no 17
 \by Gordon~E.~I. and Lyubetsky~V.~A.
 \paper Some applications of nonstandard analysis in the theory of Boolean valued measures
  \jour Dokl. Akad. Nauk SSSR
 \yr 1981
 \vol 256
 \issue 5
 \pages 1037--1041
 \endref

 \ref\no 18
 \by Gordon~E.~I. and Lyubetsky~V.~A.
 \paper Boolean completion of uniformities
  \inbook  {\it Studies on Nonclassical Logics and Formal Systems\/} [Russian]
   \publaddr Moscow
 \publ Nauka
 \yr 1983
 \pages 82--153
 \endref

 \ref\no 19
 \by Ando~T.,
 \paper Banachverb\"ande und positive Projektionen
\jour Math.~Z.
\vol 109
\pages 121--130
\year 1969
 \endref

 \ref\no 20
 \by Kusraev~A.~G.
 \paper Atomicity in injective Banach lattices
 \jour Vladikavkaz. Mat. Zh.
 \yr 2015
 \vol 17
 \issue 3
 \pages 34--42
 \endref
 }
\endRefs
\enddocument